\documentclass{amsart}

\usepackage{amssymb}
\usepackage{amsmath}
\usepackage{amsthm}
\usepackage{amsfonts}

\usepackage{a4wide}

\newtheorem{lemma}{Lemma}

\newtheorem{theorem}{Theorem}

\theoremstyle{definition}

\newtheorem{example}{Example}

\begin{document}


\title[First excited state]
{First excited state with \\ moderate rank distribution}

\author{Daniel C. Mayer}
\address{Naglergasse 53\\8010 Graz\\Austria}
\email{algebraic.number.theory@algebra.at}
\urladdr{http://www.algebra.at}

\thanks{Research supported by the Austrian Science Fund (FWF): projects J0497-PHY, P26008-N25, and by EUREA}



\date{Wednesday, 13 October 2021}


\begin{abstract}
Evidence is provided for the existence of
infinite periodic sequences of Schur \(\sigma\)-groups \(G\)
with commutator quotient \(G/G^\prime\simeq C_{3^e}\times C_3\), \(e\ge 7\),
and logarithmic order \(\mathrm{lo}(G)=10+e\).
With respect to their maximal subgroups \(H_1,\ldots,H_3;H_4\),
they have \textit{moderate} rank distribution \(\varrho(G)=(\mathrm{rank}_3(H_i/H_i^\prime))_{1\le i\le 4}\sim (2,2,3;3)\)
and represent the \textit{first excited state} of their punctured transfer kernel types \(\varkappa(G)\),
which is characterized by
a polarized component of the abelian quotient invariants \(\alpha_1(G)=(H_i/H_i^\prime)_{1\le i\le 4}\)
with \(\mathrm{lo}=6+e\) in contrast to the ground state with \(\mathrm{lo}=4+e\).
\end{abstract}

\maketitle


\section{Introduction}
\label{s:Intro}

\noindent
This is the third (and last) of a series of three articles
devoted to periodic sequences of Schur \(\sigma\)-groups \(G\)
\cite{KoVe1975,Ag1998,BuMa2015,BBH2017}
with bicyclic commutator quotients \(G/G^\prime\simeq C_{3^e}\times C_3\)
having one non-elementary component with logarithmic exponent \(e\ge 2\).
The periodicity appears in the shape of an infinite chain of 
immediate \(p\)-descendants of finite \(3\)-groups
with variable \(e\ge e_0\) bigger than a starting value.
The Schur \(\sigma\)-groups arise as leaves
of finite twigs with constant structure
emanating from the vertices of the chain.

In the first article
\cite{Ma2021a}
of the trilogy,
periodicity of pairs of \textit{metabelian} Schur \(\sigma\)-groups
sets in with \(e_0=3\),
and the \textit{ground state} of non-metabelian Schur \(\sigma\)-groups \(G\)
with moderate rank distribution
\(\varrho(G)\in\lbrace (2,2,2;3),(2,2,3;3)\rbrace\)
begins to become periodic for \(e_0=5\).
The typical bifurcation between \(G\) and its metabelianization \(M=G/G^{\prime\prime}\)
\textit{degenerates} to a simple \(p\)-descendant relation, already for \(e\ge 4\), i.e.,
the siblings topology becomes a child topology
\cite{Ma2016b}.

The primary motivation for the investigations in the present article
was the question what will happen with the more complicated fork topology
between \(G\) and \(M=G/G^{\prime\prime}\)
for the \textit{first excited state} of non-metabelian Schur \(\sigma\)-groups \(G\)
with moderate rank distribution
\(\varrho(G)\in\lbrace (2,2,2;3),\) \((2,2,3;3)\rbrace\),
which was conjectured to become periodic for \(e_0=7\)
in the conclusion of
\cite[\S\ 12]{Ma2021a}.
We shall see that \(e_0=7\) is confirmed, and the bifurcation
degenerates to an iterated \(p\)-descendant relation, already for \(e\ge 6\).

The most difficult situation of non-metabelian Schur \(\sigma\)-groups \(G\)
with \textit{elevated} rank distribution
\(\varrho(G)=(3,3,3;3)\)
was completely clarified in
\cite{Ma2021b}.
Periodicity sets in with \(e_0=9\) and the extremely complicated fork topology
freezes to a \textit{common bifurcation of infinite order} for all values \(e\ge 4\).
An additional complication with decisive negative impact
on experimental arithmetical realizations by \(3\)-class tower groups
\(G\simeq\mathrm{Gal}(\mathrm{F}_3^\infty(K)/K)\)
of imaginary quadratic number fields \(K\)
is the requirement of logarithmic abelian quotient invariants \(\alpha_2(G)\) of second order,
whereas in the present article and in
\cite{Ma2021a}
we have the immense benefit that invariants \(\alpha_1(G)\) of first order suffice. 


\section{Main Theorems: periodic Schur \(\sigma\)-groups}
\label{s:Main}

\noindent
In order to emphasize logical independence,
we partition our main result into
existence, uniqueness, explicit construction in virtue of periodicity,
and structural invariants.

\begin{theorem}
\label{thm:Existence}
\textbf{(Existence Theorem.)}
For each logarithmic exponent \(e\ge 2\),
and for each of three punctured transfer kernel types \((\mathrm{pTKT})\),
\begin{equation}
\label{eqn:pTKT}
\mathrm{D}.5,\ \varkappa(G)\sim (112;3), \quad
\mathrm{C}.4,\ \varkappa(G)\sim (113;3), \quad
\mathrm{D}.10,\ \varkappa(G)\sim (114;3),
\end{equation}
there exists a unique pair of Schur \(\sigma\)-groups \(G\)
with commutator quotient \(G/G^\prime\simeq C_{3^e}\times C_3\)
and logarithmic order \(\mathrm{lo}(G)=10+e\).
\end{theorem}

\noindent
Observe that existence is warranted also
in the pre-periodic range \(2\le e\le 6\).

\begin{theorem}
\label{thm:Periodicity}
\textbf{(Periodicity Theorem.)}
For each logarithmic exponent \(e\ge 7\),
the unique pair of Schur \(\sigma\)-groups \(G\)
with commutator quotient \(G/G^\prime\simeq C_{3^e}\times C_3\)
and logarithmic order \(\mathrm{lo}(G)=10+e\)
is given explicitly by the periodic sequence of descendants
(notation according to
\cite{BEO2005,GNO2006})
\begin{equation}
\label{eqn:Periodicity}
G\simeq W_\ell\lbrack-\#1;1\rbrack^{e-7}-\#1;i-\#1;1-\#1;1, \quad i\in\lbrace 2,3\rbrace
\end{equation}
with periodic root \(W_\ell=\mathrm{SmallGroup}(6561,93)-\#2;1-\#2;1-\#2;\ell\), where
\begin{equation}
\label{eqn:Parameter}
\begin{aligned}
\ell &= 2 \text{ for type }
\mathrm{D}.10,\ \varkappa(G)\sim (114;3), \\
\ell &= 4 \text{ for type }
\mathrm{C}.4,\ \varkappa(G)\sim (113;3), \\
\ell &= 5 \text{ for type }
\mathrm{D}.5,\ \varkappa(G)\sim (112;3).
\end{aligned}
\end{equation}
The metabelianization of \(G\) is given by
\begin{equation}
\label{eqn:Metabelianization}
M=G/G^{\prime\prime}\simeq W_\ell\lbrack-\#1;1\rbrack^{e-7}-\#1;i, \quad i\in\lbrace 2,3\rbrace,
\end{equation}
of log order \(\mathrm{lo}(G)=8+e\), in fact,
\(G^{\prime\prime}\) is cyclic of order \(9\) and is contained in the center of \(G\).
\end{theorem}

\begin{theorem}
\label{thm:Structure}
\textbf{(Structure Theorem.)}
For each logarithmic exponent \(e\ge 4\),
the unique pair of Schur \(\sigma\)-groups \(G\)
with commutator quotient \(G/G^\prime\simeq C_{3^e}\times C_3\)
and logarithmic order \(\mathrm{lo}(G)=10+e\)
possesses logarithmic abelian quotient invariants of first order
\begin{equation}
\label{eqn:AQI1}
\begin{aligned}
\alpha_1(G) &=  \lbrack (e+1)1,(e+1)1,e33;e11\rbrack \text{ for type }
\mathrm{D}.10, \text{ and } \\
\alpha_1(G) &=  \lbrack (e+1)1,(e+1)1,(e+1)32;e11\rbrack \text{ for type }
\mathrm{C}.4 \text{ and } \mathrm{D}.5.
\end{aligned}
\end{equation}
thus, in both cases, moderate rank distribution
\(\varrho(G)\sim (2,2,3;3)\),
and soluble length \(\mathrm{sl}(G)=3\).
\end{theorem}

\noindent
Observe that the periodic invariants partially also occur
in the pre-periodic range \(4\le e\le 6\),
but in different form for \(2\le e\le 3\).
The proof will be developed in \S\
\ref{s:Proof},
illustrated by Figure
\ref{fig:Tree}.


\section{Pre-periodic Schur \(\sigma\)-groups}
\label{s:Preperiodic}

\noindent
Since periodicity in Theorem
\ref{thm:Periodicity}
sets in with \(e_0=7\),
we must provide supplementary results
for the pre-periodic range \(2\le e\le 6\).
They are also justified in \S\
\ref{s:Proof}
and illustrated by Figure
\ref{fig:Tree}.

\begin{theorem}
\label{thm:21}
The unique pair of Schur \(\sigma\)-groups \(G\)
with commutator quotient \(G/G^\prime\simeq C_9\times C_3\)
and logarithmic order \(\mathrm{lo}(G)=12\)
has first abelian quotient invariants \(\alpha_1(G)\sim (31,31,431;211)\)
and is given by
\(G\simeq\mathrm{SmallGroup}(2187,168)-\#2;7-\#1;4-\#2;i\)
with metabelianization
\(M\simeq\mathrm{SmallGroup}(2187,168)-\#1;7-\#1;4-\#1;i\),
where \(i\in\lbrace 2,9\rbrace\) for type \(\mathrm{D}.5\),
\(i\in\lbrace 3,8\rbrace\) for type \(\mathrm{C}.4\), and
\(i\in\lbrace 5,6\rbrace\) for type \(\mathrm{D}.10\).
Thus, \(F=\mathrm{SmallGroup}(2187,168)\) is fork between \(M\) and \(G\).
\end{theorem}


\begin{theorem}
\label{thm:31}
The unique pair of Schur \(\sigma\)-groups \(G\)
with commutator quotient \(G/G^\prime\simeq C_{27}\times C_3\)
and logarithmic order \(\mathrm{lo}(G)=13\)
has first abelian quotient invariants \(\alpha_1(G)\sim (41,41,432;311)\)
and is given by
\(G\simeq\mathrm{SmallGroup}(6561,98)-\#2;1-\#1;1-\#2;i\)
with metabelianization
\(M\simeq\mathrm{SmallGroup}(6561,98)-\#1;3-\#1;1-\#1;i\),
where \(i\in\lbrace 2,3\rbrace\) for type \(\mathrm{D}.10\),
\(i\in\lbrace 5,9\rbrace\) for type \(\mathrm{C}.4\), and
\(i\in\lbrace 6,8\rbrace\) for type \(\mathrm{D}.5\).
Thus, \(F=\mathrm{SmallGroup}(6561,98)\) is the fork between \(M\) and \(G\).
\end{theorem}


\begin{theorem}
\label{thm:41}
The unique pair of Schur \(\sigma\)-groups \(G\)
with commutator quotient \(G/G^\prime\simeq C_{81}\times C_3\)
and logarithmic order \(\mathrm{lo}(G)=14\)
has first abelian quotient invariants \(\alpha_1(G)\sim (51,51,532;411)\)
or \(\alpha_1(G)\sim (51,51,433;411)\)
and is given by
\(G\simeq F-\#2;1-\#2;i\)
with metabelianization
\(M\simeq F-\#1;2-\#1;i\),
in terms of the fork \(F=\mathrm{SmallGroup}(6561,93)-\#2;6\),
where \(i\in\lbrace 2,3\rbrace\) for type \(\mathrm{D}.10\),
\(i\in\lbrace 5,9\rbrace\) for type \(\mathrm{C}.4\), and
\(i\in\lbrace 6,8\rbrace\) for type \(\mathrm{D}.5\).
\end{theorem}


\begin{theorem}
\label{thm:51}
The unique pair of Schur \(\sigma\)-groups \(G\)
with commutator quotient \(G/G^\prime\simeq C_{243}\times C_3\)
and logarithmic order \(\mathrm{lo}(G)=15\)
has first abelian quotient invariants \(\alpha_1(G)\sim (61,61,632;511)\)
or \(\alpha_1(G)\sim (61,61,533;511)\)
and is given by
\(G\simeq F-\#2;i-\#1;1\)
with metabelianization
\(M\simeq F-\#1;(i+1)\),
in terms of the fork \(F=\mathrm{SmallGroup}(6561,93)-\#2;1-\#2;6\),
where \(i\in\lbrace 2,3\rbrace\) for type \(\mathrm{D}.10\),
\(i\in\lbrace 5,9\rbrace\) for type \(\mathrm{C}.4\), and
\(i\in\lbrace 6,8\rbrace\) for type \(\mathrm{D}.5\).
\end{theorem}


\begin{theorem}
\label{thm:61}
The unique pair of Schur \(\sigma\)-groups \(G\)
with commutator quotient \(G/G^\prime\simeq C_{729}\times C_3\)
and logarithmic order \(\mathrm{lo}(G)=16\)
has first abelian quotient invariants \(\alpha_1(G)\sim (71,71,732;611)\)
or \(\alpha_1(G)\sim (71,71,633;611)\)
and is given by
\(G\simeq M-\#1;1-\#1;1\)
as an iterated \(p\)-descendant of the metabelianization
\(M\simeq\mathrm{SmallGroup}(6561,93)-\#2;1-\#2;1-\#2;i\),
where \(i\in\lbrace 7,8\rbrace\) for type \(\mathrm{D}.10\),
\(i\in\lbrace 10,14\rbrace\) for type \(\mathrm{C}.4\), and
\(i\in\lbrace 11,13\rbrace\) for type \(\mathrm{D}.5\).
\end{theorem}


\section{Arithmetical realizations}
\label{s:Applications}

\noindent
The automorphism group \(\mathrm{Gal}(\mathrm{F}_3^\infty(K)/K)\)
of the maximal unramified pro-\(3\) extension \(\mathrm{F}_3^\infty(K)\)
of an imaginary quadratic number field \(K=\mathbb{Q}(\sqrt{d})\)
with negative fundamental discriminant \(d<0\)
must be a Schur \(\sigma\)-group
\cite{KoVe1975,BuMa2015,BBH2017}.

With the aid of the computational algebra system Magma
\cite{BCP1997,BCFS2021,MAGMA2021},
we conducted a search for fundamental discriminants \(-10^9<d<0\)
such that the \(3\)-class group \(\mathrm{Cl}_3(K)\simeq C_{3^e}\times C_3\)
is non-elementary bicyclic with \(2\le e\le 8\),
and the capitulation in the four unramified cyclic cubic extensions of \(K\)
is the first excited state of one of the three types
\(\mathrm{C}.4\), \(\mathrm{D}.5\), \(\mathrm{D}.10\)
under investigation or of the closely related type \(\mathrm{D}.6\).
The class field routines by Fieker
\cite{Fi2001}
were employed.


In order to enable comparison with analogous cases
of imaginary quadratic number field \(K\) with
elementary bicyclic \(3\)-class group \(\mathrm{Cl}_3(K)\simeq C_3\times C_3\),
we begin with a recall of arithmetical information in
\cite[Fig. 1--2, pp. 24--25]{Ma2015b}
concerning the first excited state of
capitulation types in section \(\mathrm{E}\)
which are characterized by the polarization \(43\) of \(\alpha_1(K)\).
The rank distribution is
either \(\varrho(K)\sim (2,3,2,2)\) for the former two types
or \(\varrho(K)=(2,2,2,2)\) for the latter two types.

\begin{example}
\label{exm:11}
For \(e=1\), the hits with absolutely minimal discriminants are \\
\(d=-262\,744\) for type \(\mathrm{E}.14\), \(\varkappa(K)\sim (3122)\), \(\alpha_1(K)\sim (43,111,21,21)\), \\
\(d=-268\,040\) for type \(\mathrm{E}.6\), \(\varkappa(K)\sim (1122)\), \(\alpha_1(K)\sim (43,111,21,21)\), \\
\(d=-297\,079\) for type \(\mathrm{E}.9\), \(\varkappa(K)\sim (2334)\), \(\alpha_1(K)\sim (21,43,21,21)\), \\
\(d=-370\,740\) for type \(\mathrm{E}.8\), \(\varkappa(K)\sim (2234)\), \(\alpha_1(K)\sim (21,43,21,21)\).
\end{example}


Now we come to examples for
non-elementary bicyclic \(3\)-class groups \(\mathrm{Cl}_3(K)\simeq C_{3^e}\times C_3\), \(e\ge 2\).
For \(e=2\), the polarization of \(\alpha_1(K)\) is irregular but uniform for all types.
It is given by \((e+2)31\) instead of \((e+1)32\).

\begin{example}
\label{exm:21}
For \(e=2\), the hits with absolutely minimal discriminants are \\
\(d=-210\,164\) for type \(\mathrm{D}.6\), \(\varkappa(K)\sim (123;1)\), \(\alpha_1(K)\sim (31,31,31;431)\), \\
\(d=-320\,968\) for type \(\mathrm{C}.4\), \(\varkappa(K)\sim (113;3)\), \(\alpha_1(K)\sim (31,31,431;211)\), \\
\(d=-354\,232\) for type \(\mathrm{D}.10\), \(\varkappa(K)\sim (114;3)\), \(\alpha_1(K)\sim (31,31,431;211)\), \\
\(d=-776\,747\) for type \(\mathrm{D}.5\), \(\varkappa(K)\sim (112;3)\), \(\alpha_1(K)\sim (31,31,431;211)\).
\end{example}


For \(e=3\), the polarization of \(\alpha_1(K)\) is uniform for all types.
It is given by \((e+1)32\).

\begin{example}
\label{exm:31}
For \(e=3\), the hits with absolutely minimal discriminants are \\
\(d=-642\,491\) for type \(\mathrm{D}.5\), \(\varkappa(K)\sim (112;3)\), \(\alpha_1(K)\sim (41,41,432;311)\), \\
\(d=-1\,021\,523\) for type \(\mathrm{D}.6\), \(\varkappa(K)\sim (123;1)\), \(\alpha_1(K)\sim (41,41,41;432)\), \\
\(d=-1\,052\,072\) for type \(\mathrm{C}.4\), \(\varkappa(K)\sim (113;3)\), \(\alpha_1(K)\sim (41,41,432;311)\), \\
\(d=-1\,265\,747\) for type \(\mathrm{D}.10\), \(\varkappa(K)\sim (114;3)\), \(\alpha_1(K)\sim (41,41,432;311)\).
\end{example}


For all \(e\ge 4\), the polarization of \(\alpha_1(K)\) depends on the type.
For \(\mathrm{C}.4\) and \(\mathrm{D}.5\) it is \((e+1)32\), whereas
for \(\mathrm{D}.10\) and \(\mathrm{D}.6\) we have the variant \(e33\).

\begin{example}
\label{exm:41}
For \(e=4\), the hits with absolutely minimal discriminants are \\
\(d=-2\,249\,263\) for type \(\mathrm{D}.6\), \(\varkappa(K)\sim (123;1)\), \(\alpha_1(K)\sim (51,51,51;433)\), \\
\(d=-2\,959\,235\) for type \(\mathrm{C}.4\), \(\varkappa(K)\sim (113;3)\), \(\alpha_1(K)\sim (51,51,532;411)\), \\
\(d=-4\,076\,823\) for type \(\mathrm{D}.10\), \(\varkappa(K)\sim (114;3)\), \(\alpha_1(K)\sim (51,51,433;411)\), \\
\(d=-5\,231\,284\) for type \(\mathrm{D}.5\), \(\varkappa(K)\sim (112;3)\), \(\alpha_1(K)\sim (51,51,532;411)\).
\end{example}


\begin{example}
\label{exm:51}
For \(e=5\), the hits with absolutely minimal discriminants are \\
\(d=-5\,593\,787\) for type \(\mathrm{D}.10\), \(\varkappa(K)\sim (114;3)\), \(\alpha_1(K)\sim (61,61,533;511)\), \\
\(d=-14\,885\,751\) for type \(\mathrm{C}.4\), \(\varkappa(K)\sim (113;3)\), \(\alpha_1(K)\sim (61,61,632;511)\), \\
\(d=-18\,597\,255\) for type \(\mathrm{D}.5\), \(\varkappa(K)\sim (112;3)\), \(\alpha_1(K)\sim (61,61,632;511)\), \\
\(d=-18\,731\,096\) for type \(\mathrm{D}.6\), \(\varkappa(K)\sim (123;1)\), \(\alpha_1(K)\sim (61,61,61;533)\).
\end{example}


\begin{example}
\label{exm:61}
For \(e=6\), the hits with absolutely minimal discriminants are \\
\(d=-11\,591\,183\) for type \(\mathrm{D}.5\), \(\varkappa(K)\sim (112;3)\), \(\alpha_1(K)\sim (71,71,732;611)\), \\
\(d=-17\,740\,111\) for type \(\mathrm{C}.4\), \(\varkappa(K)\sim (113;3)\), \(\alpha_1(K)\sim (71,71,732;611)\), \\
\(d=-33\,942\,367\) for type \(\mathrm{D}.6\), \(\varkappa(K)\sim (123;1)\), \(\alpha_1(K)\sim (71,71,71;633)\).
\end{example}


\begin{example}
\label{exm:71}
For \(e=7\), the hits with absolutely minimal discriminants are \\
\(d=-111\,733\,415\) for type \(\mathrm{D}.10\), \(\varkappa(K)\sim (114;3)\), \(\alpha_1(K)\sim (81,81,733;711)\), \\
\(d=-116\,407\,871\) for type \(\mathrm{D}.5\), \(\varkappa(K)\sim (112;3)\), \(\alpha_1(K)\sim (81,81,832;711)\).
\end{example}


\begin{example}
\label{exm:81}
For \(e=8\), the hits with absolutely minimal discriminants are \\
\(d=-98\,311\,919\) for type \(\mathrm{D}.5\), \(\varkappa(K)\sim (112;3)\), \(\alpha_1(K)\sim (91,91,932;811)\).
\end{example}


\begin{theorem}
\label{thm:ThreeStage}
\textbf{(Three Stage Tower Theorem.)}
Any imaginary quadratic number field \(K=\mathbb{Q}(\sqrt{d})\), \(d<0\),
with non-elementary bicyclic \(3\)-class group \(\mathrm{Cl}_3(K)\simeq C_{3^e}\times C_3\), \(e\ge 2\),
and Artin pattern \(\mathrm{AP}(K)=(\varkappa(K),\alpha_1(K))\) given by Formulas
\eqref{eqn:pTKT} for \(\varkappa(K)\) and \eqref{eqn:AQI1} for \(\alpha_1(K)\)
has a finite \(3\)-class field tower \(\mathrm{F}_3^\infty(K)\) with precisely three stages, \(\ell_3(K)=3\).
For \(e\ge 7\), the \(3\)-class tower group \(G=\mathrm{Gal}(\mathrm{F}_3^\infty(K)/K)\)
is given by Formula \eqref{eqn:Periodicity},
and the second \(3\)-class group \(M=G/G^{\prime\prime}\simeq\mathrm{Gal}(\mathrm{F}_3^2(K)/K)\)
by Formula \eqref{eqn:Metabelianization},
both in dependence on Formula \eqref{eqn:Parameter}.
\end{theorem}

\begin{proof}
According to the Galois correspondence of field theory
and the Artin reciprocity law of class field theory
\cite{Ar1927},
the maximal self-conjugate subgroups
\(H_1,\ldots,H_3;H_4\)
of the automorphism group \(G=\mathrm{Gal}(\mathrm{F}_3^\infty(K)/K)\)
of the maximal unramified pro-\(3\) extension of an algebraic number field \(K\)
with bicyclic \(3\)-class group \(\mathrm{Cl}_3(K)\simeq C_{3^e}\times C_3\)
correspond to the unramified cyclic cubic extensions
\(L_1,\ldots,L_3;L_4\) of \(K\),
and the abelian quotient invariants
\(\alpha_1(G)=\lbrack G/G^\prime;(H_i/H_i^\prime)_{1\le i\le 4}\rbrack\)
coincide with the abelian type invariants of \(3\)-class groups
\(\alpha_1(K)=\lbrack\mathrm{Cl}_3(K);(\mathrm{Cl}_3(L_i))_{1\le i\le 4}\rbrack\).
According to Artin's theory of the transfer
\cite{Ar1929},
the Schur transfer homomorphisms
\(T_i:\,G/G^\prime\to H_i/H_i^\prime\)
correspond to the extension homomorphisms
\(\tau_i:\,\mathrm{Cl}_3(K)\to\mathrm{Cl}_3(L_i)\) of \(3\)-ideal classes,
and the punctured transfer kernel type
\(\varkappa(G)=(\ker(T_i))_{1\le i\le 4}\)
coincides with the punctured capitulation type
\(\varkappa(K)=(\ker(\tau_i))_{1\le i\le 4}\).
This was discussed in more detail in
\cite{Ma2012a}
and is the foundation of the strategy of \textit{pattern recognition via Artin transfers}
\cite{Ma2020a},
which is due to the coincidence of Artin patterns
\(\mathrm{AP}(G)=(\varkappa(G),\alpha_1(G))=(\varkappa(K),\alpha_1(K))=\mathrm{AP}(K)\).
Finally, the soluble length \(\mathrm{sl}(G)\) is equal to the
length \(\ell_3(K)\) of the \(3\)-class field tower of \(K\).
For an imaginary quadratic field \(K\),
the \(3\)-class tower group \(G\) must be a Schur \(\sigma\)-group
\cite{Sh1964,KoVe1975}.
\end{proof}


\section{Proof and tree diagram}
\label{s:Proof}

\noindent
The proof of the preperiodic Theorems
\ref{thm:21}
--
\ref{thm:61}
and finally of the periodic Main Theorems
\ref{thm:Existence}
--
\ref{thm:Structure}
can be developed in accordance with the tree diagram in Figure
\ref{fig:Tree}.
All directed edges of this tree lead from descendants \(D\) to \(p\)-parents \(\pi_p(D)=D/P_{c_p-1}(D)\), \(c_p=\mathrm{cl}_p(D)\),
rather than to parents \(\pi(D)=D/\gamma_c(D)\), \(c=\mathrm{cl}(D)\).
Consequently, the figure admits actual descendant construction.

\newpage

\begin{figure}[ht]
\caption{Schur \(\sigma\)-groups \(G\) with \(\varrho(G)\sim (2,2,3;3)\), \(G/G^\prime\simeq (3^e,3)\), \(2\le e\le 9\)}
\label{fig:Tree}

{\tiny

\setlength{\unitlength}{0.9cm}
\begin{picture}(14,18.5)(-9.5,-16.3)

\put(-11,2.5){\makebox(0,0)[cb]{order}}

\put(-11,2){\line(0,-1){17}}
\multiput(-11.1,2)(0,-1){18}{\line(1,0){0.2}}


\put(-10.8,2){\makebox(0,0)[lc]{\(3^2\)}}
\put(-10.8,1){\makebox(0,0)[lc]{\(3^3\)}}
\put(-10.8,0){\makebox(0,0)[lc]{\(3^4\)}}
\put(-10.8,-1){\makebox(0,0)[lc]{\(3^5\)}}
\put(-10.8,-2){\makebox(0,0)[lc]{\(3^6\)}}
\put(-10.8,-3){\makebox(0,0)[lc]{\(3^7\)}}
\put(-10.8,-4){\makebox(0,0)[lc]{\(3^8\)}}
\put(-10.8,-5){\makebox(0,0)[lc]{\(3^9\)}}
\put(-10.8,-6){\makebox(0,0)[lc]{\(3^{10}\)}}
\put(-10.8,-7){\makebox(0,0)[lc]{\(3^{11}\)}}
\put(-10.8,-8){\makebox(0,0)[lc]{\(3^{12}\)}}
\put(-10.8,-9){\makebox(0,0)[lc]{\(3^{13}\)}}
\put(-10.8,-10){\makebox(0,0)[lc]{\(3^{14}\)}}
\put(-10.8,-11){\makebox(0,0)[lc]{\(3^{15}\)}}
\put(-10.8,-12){\makebox(0,0)[lc]{\(3^{16}\)}}
\put(-10.8,-13){\makebox(0,0)[lc]{\(3^{17}\)}}
\put(-10.8,-14){\makebox(0,0)[lc]{\(3^{18}\)}}
\put(-10.8,-15){\makebox(0,0)[lc]{\(3^{19}\)}}


\put(-11,-15){\vector(0,-1){1.5}}

\put(-9,2){\circle{0.2}}
\put(-9,2){\circle*{0.1}}

\put(-9,1){\circle{0.2}}
\put(-9,-1){\circle{0.2}}
\put(-9,-2){\circle*{0.2}}
\put(-9,-3){\circle{0.2}}
\put(-9,-4){\circle{0.2}}
\put(-9,-5){\circle{0.2}}

\put(-7,0){\circle{0.2}}
\put(-7,-2){\circle{0.2}}
\put(-7,-3){\circle*{0.2}}
\put(-7,-4){\circle{0.2}}
\put(-7,-5){\circle{0.2}}
\put(-7,-6){\circle{0.2}}

\put(-5,-2){\circle{0.2}}
\put(-5,-4){\circle*{0.2}}
\put(-5,-5){\circle{0.2}}
\put(-5,-6){\circle{0.2}}
\put(-5,-7){\circle{0.2}}

\put(-3,-4){\circle{0.2}}
\put(-3,-6){\circle*{0.2}}
\put(-3,-7){\circle{0.2}}
\put(-3,-8){\circle{0.2}}

\put(-1,-6){\circle{0.2}}
\put(-1,-8){\circle*{0.2}}
\put(-1,-9){\circle{0.2}}

\put(1,-8){\circle{0.2}}
\put(1,-10){\circle{0.2}}

\put(3,-10){\circle{0.2}}
\put(3,-11){\circle{0.2}}

\put(4,-11){\circle{0.2}}
\put(4,-12){\circle{0.2}}

\put(5,-12){\circle{0.2}}
\put(5,-13){\circle{0.2}}

\put(-9.6,-4.1){\framebox(0.2,0.2){}}
\put(-9.6,-5.1){\framebox(0.2,0.2){}}
\put(-9.5,-5){\circle*{0.1}}
\put(-9.6,-7.1){\framebox(0.2,0.2){}}
\put(-9.5,-7){\circle{0.1}}

\put(-7.6,-5.1){\framebox(0.2,0.2){}}
\put(-7.6,-6.1){\framebox(0.2,0.2){}}
\put(-7.5,-6){\circle*{0.1}}
\put(-7.6,-8.1){\framebox(0.2,0.2){}}
\put(-7.5,-8){\circle{0.1}}

\put(-5.6,-6.1){\framebox(0.2,0.2){}}
\put(-5.6,-7.1){\framebox(0.2,0.2){}}
\put(-5.5,-7){\circle*{0.1}}
\put(-5.6,-9.1){\framebox(0.2,0.2){}}
\put(-5.5,-9){\circle{0.1}}

\put(-3.6,-8.1){\framebox(0.2,0.2){}}
\put(-3.5,-8){\circle*{0.1}}
\put(-3.6,-10.1){\framebox(0.2,0.2){}}
\put(-3.5,-10){\circle{0.1}}

\put(-1.6,-10.1){\framebox(0.2,0.2){}}
\put(-1.6,-11.1){\framebox(0.2,0.2){}}
\put(-1.5,-11){\circle{0.1}}

\put(0.9,-11.1){\framebox(0.2,0.2){}}
\put(0.9,-12.1){\framebox(0.2,0.2){}}
\put(1,-12){\circle{0.1}}

\put(2.9,-12.1){\framebox(0.2,0.2){}}
\put(2.9,-13.1){\framebox(0.2,0.2){}}
\put(3,-13){\circle{0.1}}

\put(3.9,-13.1){\framebox(0.2,0.2){}}
\put(3.9,-14.1){\framebox(0.2,0.2){}}
\put(4,-14){\circle{0.1}}

\put(4.9,-14.1){\framebox(0.2,0.2){}}
\put(4.9,-15.1){\framebox(0.2,0.2){}}
\put(5,-15){\circle{0.1}}


\put(-9,2){\line(0,-1){1}}
\put(-9,1){\line(0,-1){2}}
\put(-9,-1){\line(0,-1){1}}
\put(-9,-2){\line(0,-1){1}}
\put(-9,-3){\line(0,-1){1}}
\put(-9,-4){\line(0,-1){1}}
\put(-9,-2){\line(-1,-4){0.5}}
\put(-9.5,-4){\line(0,-1){1}}
\put(-9.5,-5){\line(0,-1){2}}

\put(-9,2){\line(1,-1){2}}
\put(-7,0){\line(0,-1){2}}
\put(-7,-2){\line(0,-1){1}}
\put(-7,-3){\line(0,-1){1}}
\put(-7,-3){\line(-1,-4){0.5}}
\put(-7,-4){\line(0,-1){1}}
\put(-7,-5){\line(0,-1){1}}
\put(-7.5,-5){\line(0,-1){1}}
\put(-7.5,-6){\line(0,-1){2}}

\put(-7,0){\line(1,-1){2}}
\put(-5,-2){\line(0,-1){2}}
\put(-5,-4){\line(0,-1){1}}
\put(-5,-5){\line(0,-1){1}}
\put(-5,-6){\line(0,-1){1}}
\put(-5,-4){\line(-1,-4){0.5}}
\put(-5.5,-6){\line(0,-1){1}}
\put(-5.5,-7){\line(0,-1){2}}

\put(-5,-2){\line(1,-1){2}}
\put(-3,-4){\line(0,-1){2}}
\put(-3,-6){\line(-1,-4){0.5}}
\put(-3,-6){\line(0,-1){1}}
\put(-3,-7){\line(0,-1){1}}
\put(-3.5,-8){\line(0,-1){2}}

\put(-3,-4){\line(1,-1){2}}
\put(-1,-6){\line(0,-1){2}}
\put(-1,-8){\line(0,-1){1}}
\put(-1,-8){\line(-1,-4){0.5}}
\put(-1.5,-10){\line(0,-1){1}}

\put(-1,-6){\line(1,-1){2}}
\put(1,-8){\line(0,-1){2}}
\put(1,-10){\line(0,-1){1}}
\put(1,-11){\line(0,-1){1}}

\put(1,-8){\line(1,-1){2}}
\put(3,-10){\line(0,-1){1}}
\put(3,-11){\line(0,-1){1}}
\put(3,-12){\line(0,-1){1}}

\put(3,-10){\line(1,-1){1}}
\put(4,-11){\line(0,-1){1}}
\put(4,-12){\line(0,-1){1}}
\put(4,-13){\line(0,-1){1}}

\put(4,-11){\line(1,-1){1}}
\put(5,-12){\line(0,-1){1}}
\put(5,-13){\line(0,-1){1}}
\put(5,-14){\line(0,-1){1}}

\put(5,-12){\line(1,-1){1}}


\put(-8.8,2){\makebox(0,0)[lc]{\(\langle 2\rangle\)}}

\put(-8.8,1){\makebox(0,0)[lc]{\(\langle 3\rangle\)}}
\put(-8.8,-1){\makebox(0,0)[lc]{\(\langle 8\rangle\)}}
\put(-8.8,-2){\makebox(0,0)[lc]{\(\langle 54\rangle\)}}
\put(-8.8,-3){\makebox(0,0)[lc]{\(1;3\)}}
\put(-8.8,-4){\makebox(0,0)[lc]{\(1;1\)}}
\put(-8.8,-5){\makebox(0,0)[lc]{\(1;i\)}}

\put(-6.8,0){\makebox(0,0)[lc]{\(\langle 3\rangle\)}}
\put(-6.8,-2){\makebox(0,0)[lc]{\(\langle 13\rangle\)}}
\put(-6.8,-3){\makebox(0,0)[lc]{\(\langle 168\rangle\)}}
\put(-6.8,-4){\makebox(0,0)[lc]{\(1;7\)}}
\put(-6.8,-5){\makebox(0,0)[lc]{\(1;4\)}}
\put(-6.8,-6){\makebox(0,0)[lc]{\(1;i\)}}

\put(-4.8,-2){\makebox(0,0)[lc]{\(\langle 7\rangle\)}}
\put(-4.8,-4){\makebox(0,0)[lc]{\(\langle 98\rangle\)}}
\put(-4.8,-5){\makebox(0,0)[lc]{\(1;3\)}}
\put(-4.8,-6){\makebox(0,0)[lc]{\(1;1\)}}
\put(-4.8,-7){\makebox(0,0)[lc]{\(1;i\)}}

\put(-2.8,-4){\makebox(0,0)[lc]{\(\langle 93\rangle\)}}
\put(-2.8,-6){\makebox(0,0)[lc]{\(2;6\)}}
\put(-2.8,-7){\makebox(0,0)[lc]{\(1;2\)}}
\put(-2.8,-8){\makebox(0,0)[lc]{\(1;i\)}}

\put(-0.8,-6){\makebox(0,0)[lc]{\(2;1\)}}
\put(-0.8,-8){\makebox(0,0)[lc]{\(2;6\)}}
\put(-0.8,-9){\makebox(0,0)[lc]{\(1;i+1\)}}

\put(1.2,-6.5){\vector(-1,-1){1.5}}
\put(1.3,-6.5){\makebox(0,0)[lc]{last semi-metabelian bifurcation}}

\put(1.2,-8){\makebox(0,0)[lc]{\(2;1\)}}
\put(1.2,-10){\makebox(0,0)[lc]{\(2;i\)}}

\put(3.2,-10){\makebox(0,0)[lc]{\(2;2\), periodic root \(W_2\)}}
\put(3.2,-11){\makebox(0,0)[lc]{\(1;i\)}}

\put(4.2,-11){\makebox(0,0)[lc]{\(1;1\)}}
\put(4.2,-12){\makebox(0,0)[lc]{\(1;i\)}}

\put(5.2,-12){\makebox(0,0)[lc]{\(1;1\)}}
\put(5.2,-13){\makebox(0,0)[lc]{\(1;i\)}}

\put(6.2,-13){\makebox(0,0)[lc]{etc.}}

\put(-9.7,-4){\makebox(0,0)[rc]{\(2;3\)}}
\put(-9.7,-5){\makebox(0,0)[rc]{\(1;1\)}}
\put(-9.7,-7){\makebox(0,0)[rc]{\(2;i\)}}
\put(-9.5,-7.5){\makebox(0,0)[cc]{\(i=4,6\)}}
\put(-9.5,-8){\makebox(0,0)[lc]{\(\mathrm{E}.9\)}}
\put(-9.5,-8.5){\makebox(0,0)[lc]{\((3,3)\)}}

\put(-7.7,-5){\makebox(0,0)[rc]{\(2;7\)}}
\put(-7.7,-6){\makebox(0,0)[rc]{\(1;4\)}}
\put(-7.7,-8){\makebox(0,0)[rc]{\(2;i\)}}
\put(-7.5,-8.5){\makebox(0,0)[cc]{\(i=5,6\)}}
\put(-7.5,-9){\makebox(0,0)[lc]{\(\mathrm{D}.10\)}}
\put(-7.5,-9.5){\makebox(0,0)[lc]{\((9,3)\)}}

\put(-5.7,-6){\makebox(0,0)[rc]{\(2;1\)}}
\put(-5.7,-7){\makebox(0,0)[rc]{\(1;1\)}}
\put(-5.7,-9){\makebox(0,0)[rc]{\(2;i\)}}
\put(-5.5,-9.5){\makebox(0,0)[cc]{\(i=2,3\)}}
\put(-5.5,-10){\makebox(0,0)[lc]{\(\mathrm{D}.10\)}}
\put(-5.5,-10.5){\makebox(0,0)[lc]{\((27,3)\)}}

\put(-3.7,-8){\makebox(0,0)[rc]{\(2;1\)}}
\put(-3.7,-10){\makebox(0,0)[rc]{\(2;i\)}}
\put(-3.5,-10.5){\makebox(0,0)[cc]{\(i=2,3\)}}
\put(-3.5,-11){\makebox(0,0)[lc]{\(\mathrm{D}.10\)}}
\put(-3.5,-11.5){\makebox(0,0)[lc]{\((81,3)\)}}

\put(-2.1,-13.5){\vector(-1,4){1.3}}
\put(-2.3,-13.6){\makebox(0,0)[ct]{last non-metabelian bifurcation}}

\put(-1.7,-10){\makebox(0,0)[rc]{\(2;i\)}}
\put(-1.7,-11){\makebox(0,0)[rc]{\(1;1\)}}
\put(-1.5,-11.5){\makebox(0,0)[cc]{\(i=2,3\)}}
\put(-1.5,-12){\makebox(0,0)[lc]{\(\mathrm{D}.10\)}}
\put(-1.5,-12.5){\makebox(0,0)[lc]{\((243,3)\)}}

\put(1.2,-11){\makebox(0,0)[lc]{\(1;1\)}}
\put(1.2,-12){\makebox(0,0)[lc]{\(1;1\)}}
\put(0.5,-12.5){\makebox(0,0)[lc]{\(i=7,8\)}}
\put(0.5,-13){\makebox(0,0)[lc]{\(\mathrm{D}.10\)}}
\put(0.5,-13.5){\makebox(0,0)[lc]{\((729,3)\)}}

\put(3.2,-12){\makebox(0,0)[lc]{\(1;1\)}}
\put(3.2,-13){\makebox(0,0)[lc]{\(1;1\)}}
\put(2.5,-13.5){\makebox(0,0)[lc]{\(i=2,3\)}}
\put(2.5,-14){\makebox(0,0)[lc]{\(\mathrm{D}.10\)}}
\put(2.5,-14.5){\makebox(0,0)[lc]{\((2187,3)\)}}

\put(4.2,-13){\makebox(0,0)[lc]{\(1;1\)}}
\put(4.2,-14){\makebox(0,0)[lc]{\(1;1\)}}
\put(3.8,-14.5){\makebox(0,0)[lc]{\(i=2,3\)}}
\put(3.5,-15){\makebox(0,0)[lc]{\(\mathrm{D}.10\)}}
\put(3.5,-15.5){\makebox(0,0)[lc]{\((6561,3)\)}}

\put(5.2,-14){\makebox(0,0)[lc]{\(1;1\)}}
\put(5.2,-15){\makebox(0,0)[lc]{\(1;1\)}}
\put(4.8,-15.5){\makebox(0,0)[lc]{\(i=2,3\)}}
\put(4.5,-16){\makebox(0,0)[lc]{\(\mathrm{D}.10\)}}
\put(4.5,-16.5){\makebox(0,0)[lc]{\((19683,3)\)}}

\put(-4,1.4){\makebox(0,0)[lc]{Legend:}}

\put(-2.5,1.4){\circle{0.2}}
\put(-2.5,1.4){\circle*{0.1}}
\put(-2,1.4){\makebox(0,0)[lc]{\(\ldots\) abelian}}

\put(-2.5,1){\circle{0.2}}
\put(-2,1){\makebox(0,0)[lc]{\(\ldots\) metabelian}}

\put(-2.5,0.6){\circle*{0.2}}
\put(-2,0.6){\makebox(0,0)[lc]{\(\ldots\) metabelian with bifurcation}}

\put(-2.6,0.1){\framebox(0.2,0.2){}}
\put(-2,0.2){\makebox(0,0)[lc]{\(\ldots\) non-metabelian}}

\put(-2.6,-0.3){\framebox(0.2,0.2){}}
\put(-2.5,-0.2){\circle*{0.1}}
\put(-2,-0.2){\makebox(0,0)[lc]{\(\ldots\) non-metabelian with bifurcation}}

\put(-2.6,-0.7){\framebox(0.2,0.2){}}
\put(-2.5,-0.6){\circle{0.1}}
\put(-2,-0.6){\makebox(0,0)[lc]{\(\ldots\) non-metabelian Schur \(\sigma\)}}

\end{picture}

}

\end{figure}


\noindent
The \(p\)-group generation algorithm
\cite{HEO2005}
by Newman
\cite{Nm1977}
and O'Brien
\cite{Ob1990}
is implemented in the ANUPQ package
\cite{GNO2006}
of the computational algebra system Magma
\cite{MAGMA2021,BCFS2021,BCP1997}.
This algorithm is used to construct
all immediate \(p\)-descendants of an assigned finite \(p\)-group.
Repeated recursive applications of the algorithm,
guided by the strategy of pattern recognition via Artin transfers
\cite{Ma2020a},
eventually produce Figure
\ref{fig:Tree}.
In each step, only \(\sigma\)-descendants are allowed to pass the filter.
The figure shows an \textit{infinite main trunk}
and \textit{finite twigs} emanating from the vertices of the trunk.
Propagation along the trunk is exo-genetic with increasing commutator quotient
\((3^e,3)\mapsto (3^{e+1},3)\),
whereas all vertices of a twig share a common abelianization
and the propagation is endo-genetic.
The leftmost twig is included in order to point out analogy to the
\textit{elementary} commutator quotient \((3,3)\).
It was computed in
\cite[Fig. 6, p. 110]{Ma2018},
where historical information is provided for type \(\mathrm{E}.9\), \(\varkappa\sim (2231)\sim (3231)\), in
\cite[\S\ 4, pp. 107--111]{Ma2018}.

\newpage

\noindent
All the other twigs concern \textit{non-elementary} commutator quotients \((3^e,3)\), \(2\le e\le 9\),
restricted to the particular punctured transfer kernel type \(\mathrm{D}.10\), \(\varkappa\sim (114;3)\).
Figure
\ref{fig:Tree}
remains unchanged for the other two types \(\mathrm{C}.4\) and \(\mathrm{D}.5\)
when the parameter \(i\) is selected according to the preperiodic Theorems
\ref{thm:21}
--
\ref{thm:61}
for \(2\le e\le 6\),
and the \textit{periodic root} \(W_\ell\) for \(e\ge 7\) is replaced according to Formula
\eqref{eqn:Parameter}
in the periodic Main Theorem
\ref{thm:Periodicity}.
More changes are required for type \(\mathrm{D}.6\), \(\varkappa\sim (123;1)\),
which is only included in the number theoretic \S\
\ref{s:Applications}
but not in the group theoretic \S\S\
\ref{s:Main}
and
\ref{s:Preperiodic}.
For instance, the coclass tree with root \(\langle 729,13\rangle\) and bifucation at \(\langle 2187,168\rangle\)
\cite[Fig. 5]{Ma2021a},
which is responsible for three types \(\mathrm{D}.10\), \(\mathrm{C}.4\) and \(\mathrm{D}.5\),
must be replaced by the coclass tree with root \(\langle 729,21\rangle\) and bifucation at \(\langle 2187,191\rangle\)
\cite[Fig. 6]{Ma2021a}
or with root \(\langle 729,18\rangle\) and bifucation at \(\langle 2187,181\rangle\).
Both of the latter coclass trees give rise to the single relevant type \(\mathrm{D}.6\).
The initialization of the construction process at \(e=2\) is described in
\cite[\S\ 5]{Ma2021b},
but now the \textit{scaffold type} \(\mathrm{b}.31\), \(\varkappa\sim (044;4)\),
must be replaced by type \(\mathrm{d}.10\), \(\varkappa\sim (110;3)\), associated with
types \(\mathrm{D}.10\), \(\mathrm{C}.4\) and \(\mathrm{D}.5\).
The search immediately leads to the root \(\langle 729,13\rangle\) and the fork \(\langle 2187,168\rangle\).
Capable descendants of both step sizes \(s\in\lbrace 1,2\rbrace\)
have relative identifiers \(\lbrace 1,4,7\rbrace\),
but only \(7\) is relevant for the desired types.

In the leftmost three twigs with \(1\le e\le 3\),
the fork topologies
\cite{Ma2016b},
which begin at the bifurcation,
are isomorphic as directed graphs.
However, in the next two twigs with \(4\le e\le 5\),
the fork topologies shrink gradually,
and for all \(e\ge 6\),
the bifurcation vanishes
leaving a descendant topology.
The vertices on twigs are BCF-groups
\cite{Ma2021b}.

The vertices of the main trunk are CF-groups
\cite{AHL1977}
with type \(\mathrm{a}.1\), \(\varkappa=(000;0)\), and \(\varrho\sim (2,2,3;3)\).
For \(2\le e\le 5\), the first vertex of the twig has scaffold type \(\mathrm{d}.10\), \(\varkappa\sim (110;3)\),
but for \(e\ge 6\), the twigs entirely consist of vertices sharing a common type,
\(\mathrm{D}.10\), \(\mathrm{C}.4\) or \(\mathrm{D}.5\).

Up to \(e\le 6\),
the construction process is straightforward,
but for \(e=7\) considerable difficulties arise,
because it is hard to determine the next vertex on the main trunk.

Since an attempt with hypothetical next main trunk root
\(\langle 6561,93\rangle-\#2;1-\#2;1-\#2;1\)
and \(1709\) descendants \(D\) up to \(p\)-class \(\mathrm{cl}_p(D)\le 11\) only led to
three pairs of Schur \(\sigma\)-groups with commutator quotient \((2187,3)\)
and \(\mathrm{lo}=20\) in the \textit{second} excited state
\(\alpha_1\sim (81,81,744;711)\) for type \(\mathrm{D}.10\), respectively
\(\alpha_1\sim (81,81,843;711)\) for type \(\mathrm{C}.4\) and \(\mathrm{D}.5\),
we returned to a tour de force computation starting at the previous vertex
\(\langle 6561,93\rangle-\#2;1-\#2;1\)
on the main trunk.


\begin{lemma}
\label{lem:BruteForce}
Among the \(1708\) descendants \(D\) up to \(p\)-class \(\mathrm{cl}_p(D)\le 10\)
of the main trunk vertex \(\langle 6561,93\rangle-\#2;1-\#2;1\),
there occur the following Schur \(\sigma\)-groups \(S\):
\begin{enumerate}
\item
the expected three pairs with \(S/S^\prime\simeq (729,3)\) and \(\mathrm{lo}(S)=16\)
in the first excited state
\(\alpha_1\sim (71,71,633;611)\) for type \(\mathrm{D}.10\), respectively
\(\alpha_1\sim (71,71,732;611)\) for type \(\mathrm{C}.4\), \(\mathrm{D}.5\),
\item
the desired three pairs with \(S/S^\prime\simeq (2187,3)\) and \(\mathrm{lo}(S)=17\)
in the first excited state
\(\alpha_1\sim (81,81,733;711)\) for type \(\mathrm{D}.10\), respectively
\(\alpha_1\sim (81,81,832;711)\) for type \(\mathrm{C}.4\), \(\mathrm{D}.5\),
\item
three unexpected pairs with \(S/S^\prime\simeq (729,3)\) and \(\mathrm{lo}(S)=19\)
in the second excited state
\(\alpha_1\sim (71,71,644;611)\) for type \(\mathrm{D}.10\), respectively
\(\alpha_1\sim (71,71,743;611)\) for type \(\mathrm{C}.4\), \(\mathrm{D}.5\),
\item
and (as a superfluous byproduct) a quartet of type \(\mathrm{B}.2\), \(\varkappa(S)\sim (111;2)\),
with \(S/S^\prime\simeq (729,3)\) and \(\mathrm{lo}(S)=19\) in the first excited state \(\alpha_1\sim (71,71,732;611)\).
\end{enumerate}
\end{lemma}

\noindent
A back track search, starting from item (2) in Lemma
\ref{lem:BruteForce},
finally reveals three suitable periodic roots \(W_\ell\)
given in Formula
\eqref{eqn:Parameter}
of Theorem
\ref{thm:Periodicity}.


\section{Conclusion}
\label{s:Conclusion}

\noindent
In the present article,
our hypothesis in
\cite[\S\ 12]{Ma2021a}
that periodicity of Schur \(\sigma\)-groups \(G\)
with \(G/G^\prime\simeq C_{3^e}\times C_3\)
and one of the punctured transfer kernel types \(\mathrm{D}.10\), \(\mathrm{C}.4\) and \(\mathrm{D}.5\)
in the \textit{first excited state}
will set in with \(e\ge e_0=7\)
was verified.

The question concerning the evolution of the fork topology
between \(G\) and its metabelianization \(M=G/G^{\prime\prime}\),
for which we were really unable to make any prediction,
was answered by a gradual shrinking of the twigs in Figure
\ref{fig:Tree}
for \(4\le e\le 7\), coming along with
a last non-metabelian (second) bifurcation for \(e=4\)
and a last semi-metabelian (first) bifurcation for \(e=5\).
The fork topology completely degenerates to a
descendant topology for all \(e\ge 6\)
and periodicity of the shortest possible twigs sets in with \(e\ge e_0=7\).


\section{Acknowledgements}
\label{s:Gratifications}

\noindent
The author acknowledges
support by the Austrian Science Fund (FWF): Projects J0497-PHY and P26008-N25,
and by the Research Executive Agency of the European Union (EUREA).



\end{document}